\documentclass[twoside]{article}

\usepackage{amssymb}
\usepackage{amsmath}

 \numberwithin{equation}{section}

\newtheorem{theorem}{Theorem}

\newtheorem{lemma}{Lemma}

\newtheorem{remark}{Remark}

\begin{document}

\title{
Boundary feedback stabilization 
by piecewise constant time delay for the wave equation
}
\author{Martin Gugat\\
Lehrstuhl 2 f\"ur Angewandte Mathematik,
\\
Martensstr. 3, 91058
Erlangen, Germany
\thanks{gugat@am.uni-erlangen.de}
\and
Marius Tucsnak
\\
Institut \'Elie Cartan Nancy (Math\'ematiques)
\\
Universit\'e Henri Poincaré Nancy 1
\\
B.P. 70239, F-54506 Vandoeuvre-les-Nancy Cedex
\\
France
\thanks{Marius.Tucsnak@loria.fr}
}
\date{ }
%
%
\maketitle


\begin{abstract}
For vibrating systems, a delay in 
the application of a feedback control can destroy the stabilizing effect of
the control.
In this paper we consider
a vibrating string that
is fixed at one end and
stabilized with a boundary feedback
with delay at the other end.

We show that for certain feedback parameters
the system is exponentially
stable with constant delays
of the form $4L/c$, $8L/c$, $12L/c$ ... 
Moreover, we show that the system is  exponentially stable
with piecewise constant delays that attain
the values $4L/c$ and $8L/c$.



\end{abstract}

{Key Words: hyperbolic pde,
feedback stabilization of pdes, delay, boundary feedback, switching delay,
wave equation,  feedback with delay, time-dependent feedback-parameter,
past observation, circular string.}

AMS Subject Classification{ 49K20,  49K25, 34H05}



\section{Introduction}
Datko et al. have  described the following
problem in the application of feedback laws:
{\em
Some second-order vibrating systems cannot tolerate small time delays in their damping
}
(see \cite{da:somes}).
In  other words:
Delays can destabilize a system that is asymptotically stable
in the absence of delays
(see \cite{da:anexa}).
The problem of instability caused by small constant delays
has also been considered in \cite{lo:condi},
where a systematic frequency domain treatment of this phenomenon
has been given
and examples for the instability
created by small delays have been presented.


In \cite{gu:bound} a constant delay with the value $2L/c$ has been considered
and it has been shown that with this delay,
exponential damping is possible for feedback parameters with sufficiently small 
absolute value that have an opposite sign as the parameters that 
generate exponential damping in the case without delay.

In this paper 
we show that for a constant delay that is
an integer multiple of $4L/c$,
exponential damping with feedback parameters 
of sufficiently small 
absolute value  is possible if the feedback parameters have the same sign as
the parameters that work in the case without delay.

Moreover, we consider piecewise constant delays with values
$4L/c$ or $8L/c$ and show that also for 
delays that switch between those two values in an arbitrary way,
the energy decays exponentially for certain feedback parameters.


For the problem considered in this paper
some progress has been made in \cite{gx:bound}
for the wave equation.
In \cite{gx:dynam} the related problem for the
Euler-Bernoulli beam has been considered.
In most studies of feedback stabilization of
second-order vibrating systems, no delays are considered:
In \cite{cs:thera}, a vibrating string is
considered and a feedback law is presented for which the energy
vanishes in finite time. In \cite{gu:optim}
it is shown that the result  from \cite{cs:thera} is stable in the sense
that also with moving boundaries, the energy is driven to zero in
finite time.
%
The problem of boundary control of the wave equation has
also been studied in \cite{ru:nonha}, \cite{li:exact}, \cite{kr:optim},
\cite{kr:onmom}, \cite{ai:famil}, \cite{tu:obser}
and the references therein.

%
%

%

This paper has the following structure:
In Section \ref{single} we define the considered system and 
in Section \ref{well} we show that it is well-posed.

In Section \ref{piece} we  show that the system is stable
with piecewise constant delays that attain
the values $4L/c$ and $8L/c$.
To our knowledge, this is the first example of
a system that is stabilized with a switching delay,
where the switching occurs between the two delay values.

In the last section we show that our feedback law is
stabilizing without delay and for a certain sequence of constant delays
with appropriately chosen feedback parameters of the same sign.
We show the exponential decay of the  energy in the system.



%
%

\section{The System}
\label{single}
Let a string of length $L>0$ and the corresponding wave speed $c>0$ be given.
Define the set $\Omega =  (0,\infty) \times (0,\; L)$.
Define the set of initial states
\[B= \{(y_0,y_1) \in H^1(0,L) \times L^2(0,L): y_0(0)=0\}.\]
Let a number $\iota\in \{0,1,2,...\}$ be given.
Assume that $\delta$ is a
piecewise constant function with $\delta(t)\in [2\frac{L}{c}, 4\iota\frac{L}{c}]$ for all $t\geq 0$.

For $(y_0,y_1) \in B$
we consider the system ${\bf S_1}:$
\begin{equation}
\label{initial1} v(0,x)= y_0(x),
\end{equation}
\begin{equation}
\label{initial1z}
v_t(0,x)= y_1(x),\; x\in (0,L)
\end{equation}
\begin{equation}
\label{wellengleichung1} v_{tt}(t,x)= c^2 v_{xx}(t,x),\;
(t,x)\in \Omega
\end{equation}
\begin{equation}
\label{RB0}
v(t,0) = 0,\; t>0
\end{equation}
%
\begin{equation}
\label{RBzero}
 v_x(t,L)  = 0,\; t \in (0, 4 \iota \frac{L}{c})
\end{equation}
and
\begin{equation}
\label{RB1} v_x(t,\,L)  =  \frac{f}{c} \;  v_t\left(t - \delta(t),\, L\right),\; t> 4\iota \frac{L}{c}.
\end{equation}
Here $f$ is a real number.
The quotient $ \frac{f}{c} $ appears in
the feedback law (\ref{RB1}) in
order to make the size of the stabilizing
feedback parameters $f$ independent of $c$.
We assume that the compatibility condition $y_0(0)=0$
is satisfied since it assures that the system has continuous states, as we show in Theorem \ref{solution}.


\section{Well-posedness of the system {$\bf S_1$}}
\label{well}

In this section, we study the
well-posedness of system {$\bf S_1$}
that is
(\ref{initial1})-(\ref{RB1}).

\begin{theorem}
\label{solution}
Assume that $\delta$ is a
piecewise constant function with $\delta(t)\in [2\frac{L}{c}, 4\iota\frac{L}{c}]$ for all $t\geq 0$.

Let $(y_0, y_1) \in B$  be given.
Define the function $\alpha$ recursively by
\begin{eqnarray}
\label{alphadefinition}
\alpha(x)    & = &
\left\{
\begin{array}{rrrll}
- \frac{1}{2}\; y_0( - x ) & + & \frac{1}{2c}   \int_0^{- x} y_1(s)\,ds,\; & x \in [-L, 0), \\
\frac{1}{2} \;y_0( x )     & + & \frac{1}{2c}   \int_0^{x} y_1(s)\,ds,\;  & x \in [0, L),
\end{array}
\right.
\end{eqnarray}
and for
$k\in \{1,2,... ,2( \iota-1)\}$
and
$x \in [L + 2k L, 3 L + 2 k L)$
by
\begin{equation}
\alpha'(x)=-\alpha'(x-2L)
\end{equation}
and for $k\in \{0,1,2,...\}$ and
$x \in [L + 4\iota L + 2k L, 3L + 4 \iota L  + 2 k L)$ by
\begin{equation}
\label{2.8}
\alpha'(x)=-\alpha'(x-2L)
+ f \alpha'(x - c\delta(\frac{x}{c})) - f \alpha'(x- 2L - c\delta(\frac{x}{c}))
\end{equation}
and the condition that $\alpha$ is continuous on the
interval $[-L,\infty)$.
Let
\begin{equation}
\label{charac7}
v(t,x)= \alpha( c t + x ) - \alpha( ct - x ),\;\;(t, x)\in
\Omega .
\end{equation}
For every finite interval $I\subset [-L,\infty)$ we have
$\alpha'\in L^2(I)$.
The function $v$ is continuous on $\Omega$ and $v_t$, $v_x\in L^1_{loc}(\Omega)$.
Define the family of test functions ${\cal T}$ as
\begin{eqnarray*}
{\cal T} & = & \{ \varphi \in C^2(\Omega):\mbox{\rm  There exists
a set $Q=[t_1,t_2]\times[x_1,x_2]\subset \Omega$ }
\\
& &
\mbox{\rm such that
the support of $\varphi$ is contained in the interior of $Q$}\}.
\end{eqnarray*}
The function $v$ satisfies the wave equation (\ref{wellengleichung1})  in the following weak sense:
\begin{equation}
\label{weak}
\int_{\Omega} v_t(t,x) \varphi_t(t,x) \,d(t,x)  = c^2\int_{\Omega} v_x(t,x) \varphi_x(t,x) \,d(t,x)\;
\mbox{\it for all }\;\varphi \in {\cal T}.
\end{equation}
The function $v$ satisfies (\ref{initial1}) and (\ref{initial1z}) and (\ref{RB0})-(\ref{RB1}).
In this sense, $v$ is the solution of the  system {$\bf S_1$}
that is
(\ref{initial1})-(\ref{RB1}).
%
\end{theorem}

{\bf Proof.}
Since $y_0'\in L^2(0,L)$,  the Sobolev imbedding Theorem implies that
$y_0$ is continuous.
Moreover, $y_1$ is in $L^2(0,L)$, thus $\alpha$ is well defined.
Now we discuss the regularity of $\alpha$.
On the intervals $[-L,0)$, $[0,L)$  and
 $[L + 2k L, 3 L + 2 k L)$ ($k\in \{0,1,2,3,...\}$
the function $\alpha$ is continuous.
Due to the definition of the set  $B$ we have
\begin{eqnarray*}
\lim_{x\rightarrow 0 -} \alpha(x) & = &
- (1/2) y_0(0)=0=
(1/2)y_0(0)= \lim_{x\rightarrow 0 + } \alpha(x), \\
\lim_{x\rightarrow L -} \alpha(x) & = &
\frac{1}{2} y_0(L) + \frac{1}{2c} \int_0^L y_1(s)\, ds
\\
&  = &   \frac{1}{c} \int_0^L y_1(s)\, ds
-
\left(
-  \frac{1}{2}  y_0(L)
+ \frac{1}{2c} \int_0^L y_1(s)\, ds \right)
\\
& = &
 \frac{1}{c} \int_0^L y_1(s)\, ds   - \alpha(-L)
=
 \lim_{x\rightarrow L +} \alpha(x), \\
\alpha(3L + 2kL) & = &
(f-1) \alpha(L + 2kL) - f \alpha(2kL - L) + C_k
\\
&  = &  \lim_{x \rightarrow  3 L + 2kL -} \alpha(x)
\end{eqnarray*}
hence $\alpha$ is continuous on the interval $[-L,\infty)$.
The derivative $\alpha'$ in the sense of distributions exists on
the intervals $(-L,0)$,  $(0,L)$, $(L, 3L)$  and $
(3L + 2k L, 5 L + 2 k L)$
as  $L^2$-function.
Since $\alpha$ is continuous, this implies that $\alpha$
is absolutely continuous on $(-L, \infty)$.
Hence $\alpha' \in L^2_{loc}(-L, \infty)$.
The continuity of $v$ follows from the continuity of $\alpha$.
For $t=0$ and  $x \in (0,L)$ we have
\[v(0,x) = \alpha(x) - \alpha(-x)= y_0(x).\]
For  $(t,x) \in \Omega$ almost everywhere, we have
\begin{equation}
\label{vtgleichung}
v_t(t,x) = c [\alpha'(x+ct ) - \alpha'(-x + c t )].
\end{equation}
Thus the definition of $\alpha$ implies the equation $v_t(0,x)= y_1(x)$.
Hence the initial conditions (\ref{initial1}) and (\ref{initial1z}) are valid.

For  $(t,x) \in \Omega$ almost everywhere, we have
\begin{equation}
\label{vxgleichung}
v_{x}(t,x)  =    \alpha'(x+ct ) + \alpha'(-x + c t ).
\end{equation}
By Tonelli's Theorem (see e.g.\cite{pe:analy}), (\ref{vxgleichung}) implies $v_x\in L^1_{loc}(\Omega)$
and  (\ref{vtgleichung}) implies $v_t\in L^1_{loc}(\Omega)$.

For all $\varphi \in {\cal T}$, integration by parts,
(\ref{vxgleichung})
and (\ref{vtgleichung})
 yield
\begin{eqnarray*}
& & \int_{\Omega} v_x(t,x) \varphi_x(t,x) \,d(t,x)
\\
& = &
\int_{x_1}^{x_2} \int_{t_1}^{t_2}
\varphi_x(t,x)  [\alpha'(x+ct ) + \alpha'(-x + c t )] \,dt\,dx
\\
& = &
- \int_{x_1}^{x_2} \int_{t_1}^{t_2}
\varphi_{xt}(t,x)  [\alpha(x+ct ) + \alpha(-x + c t )]/c \,dt\,dx
\\
& = &
-  \int_{t_1}^{t_2} \int_{x_1}^{x_2}
\varphi_{tx}(t,x)  [\alpha(x+ct ) + \alpha(-x + c t )]/c \,dx\,dt
\\
& = &
  \int_{t_1}^{t_2} \int_{x_1}^{x_2}
\varphi_{t}(t,x)  [\alpha'(x+ct ) - \alpha'(-x + c t )]/c \,dx\,dt
\\
& = &
  \int_{\Omega}
\varphi_{t}(t,x) \, v_t(t,x)/c^2 \,d(x,t)
\end{eqnarray*}
hence (\ref{weak}) holds.


For $x=0$ we have
$v(t,0)= \alpha(ct) - \alpha(ct) =0$, hence
at $x=0$ the boundary condition $v(t,0)=0$
holds for all $ t>0$.

For $x=L$, (\ref{vxgleichung}) implies
for $t \in (0, 4 \iota \frac{L}{c})$ the equation
$$
v_x(t,L)=\alpha'(L + c t )  + \alpha'( c t -L)
= - \alpha'(ct - L) + \alpha'(ct -L) =0.
$$
Therefore, the boundary condition (\ref{RBzero})
holds for all $t \in (0, 4\iota \frac{L}{c})$.


For $t > 4 \iota L/c$, we have
\begin{eqnarray*}
v_x(t,L) & = & \alpha'(c t + L) + \alpha'(c t -L)
\\
& = & f \alpha'(ct -c\delta(t)) - f \alpha'(c t - 2L
-c\delta(t))
\\
& = &
f[\alpha'(L + c t - c\delta(t)) - \alpha'(-L + c t - c\delta(t))]
 \\
& = & (f/c)\; v_t(t- c \delta(t),\,L).
\end{eqnarray*}
Therefore, the boundary condition (\ref{RB1})
holds for all $t > 4\iota \frac{L}{c}$.

\begin{remark}
Note that our system has a continuous state.
Optimal boundary control problems for the
wave equation with countinuous states
have been considered in
\cite{gu:optimc}.
The proof of Theorem \ref{solution}
is similar to the proof of
Theorem 4.1 in \cite{gu:optime}.
Theorem \ref{solution} is a generalization of
Theorem 1 in \cite{gu:bound}, where the case
$\delta(x)=2L/c$ has been considered.
\end{remark}

\subsection{Transformation of the recursion to a vector recursion}

Instead of the recursion (\ref{2.8})
we can also use the following  linear system to characterize the solution
of {$\bf S_1$}:
\begin{equation}
\label{21}
\left(
\begin{array}{r}
\alpha'(x)
\\
\alpha'(x-2L)
\\
\alpha'(x-4L)
\\
\alpha'(x-6L)
\\
\alpha'(x -c\delta(\frac{x}{c}))
\end{array}
\right)
=
\left(
\begin{array}{rrrrr}
-1 & 0 & 0 & f & -f
\\
1 & 0 & 0 & 0 & 0
\\
0 & 1 & 0 & 0 & 0
\\
0 & 0 & 1 & 0 & 0
\\
0 & 0 & 0 & 1 & 0
\end{array}
\right)
\;
\left(
\begin{array}{r}
\alpha'(x-2L)
\\
\alpha'(x-4L)
\\
\alpha'(x-6L)
\\
\alpha'(x-c\delta(\frac{x}{c}))
\\
\alpha'(x- 2L  - c\delta(\frac{x}{c}))
\end{array}
\right).
\end{equation}
Let $B_2$ be the matrix in system (\ref{21}).
Let $\det(\lambda I - B_2)=p_f(\lambda)$
denote the characteristic polynomial of $B_2$.
Then we have the equation
$$
p_f(\lambda)=\lambda^5 + \lambda^4 -\lambda f + f.$$

If $\delta(x)=4$, we can write (\ref{2.8}) in the form of the linear system
\begin{equation}
\label{20}
\left(
\begin{array}{r}
\alpha'(x)
\\
\alpha'(x-2L)
\\
\alpha'(x-4L)
\\
\alpha'(x-6L)
\\
\alpha'(x-8L)
\end{array}
\right)
=
\left(
\begin{array}{rrrrr}
-1 &  f & -f & 0 & 0
\\
1 & 0 & 0 & 0 & 0
\\
0 & 1 & 0 & 0 & 0
\\
0 & 0 & 1 & 0 & 0
\\
0 & 0 & 0 & 1 & 0
\end{array}
\right)
\;
\left(
\begin{array}{r}
\alpha'(x-2L)
\\
\alpha'(x-4L)
\\
\alpha'(x-6L)
\\
\alpha'(x-8L)
\\
\alpha'(x-10L)
\end{array}
\right).
\end{equation}
Let $B_1$ be the matrix in system (\ref{20}).
Let $\det(\lambda I - B_1)=d_f(\lambda)$
denote the characteristic polynomial of $B_1$.
Then we have the equation
$$
d_f(\lambda)=\lambda^2 (\lambda^3 + \lambda^2 -\lambda f + f).$$


\subsection{The characteristic Polynomial}
\label{charac}
Let $j$ be a natural number.
For a real number $f$ we define the polynomial

$$ p_f(t) = t^{2j+1} + t^{2j} - f t + f.$$

So for $f=0$ we have
$p_0(t)=t^{2j}(1+t)$ with the roots $(-1)$ (with multiplicity 1)
and zero as the second root (with multiplicity $2j$).

\begin{lemma}
\label{lemmaone}
There exists a number $\delta_j>0$, such that for all
$f\in (-\delta_j,0)$, all roots of $p_f$ have
a modulus that is strictly less than one.
\end{lemma}

For the proof of Lemma \ref{lemmaone}, we use an
intermediate result.
The following Lemma \ref{lemma0p} states that
outside a neigbourhood of $(-1)$, all roots of $p_f$
have a  modulus that is strictly less than one.

\begin{lemma}
\label{lemma0p}
Let $z$ be a root of $p_f$ with
$|z-(-1)|> \frac{2 |f|}{ 1 - |f|}$.
Then $$|z|<1,$$
that is the modulus of $z$  is strictly less than one.
\end{lemma}
{\bf Proof.}
The equation $p_f(z)=0$ implies
$z^{2j}(z+1) -f(z+1) + 2f=0$.
Hence
$(z^{2j} - f)(z+1)= - 2f$ which implies the
inequality
$$|z^{2j} - f| =  \frac{2 |f|}{|z+1|}.$$
We have
\begin{eqnarray*}
|z|^{2j} & \leq  & |z^{2j} - f| + |f|\\
& = &   \frac{2 |f|}{|z+1|} + |f|\\
& < &   (1 -  |f|) + |f| = 1.
\end{eqnarray*}
Hence we have $|z| < 1$ and the assertion follows.

Note that if $f<4j+1$, for all $t\geq 1$ we have
$p_f'(t) =t^{2j-1}((2j+1) t + 2j)-f > 0$, and since $p_f(1)=2$ this implies  that $p_f$ does not have a root in
$[1,\infty)$.

Now we come to the proof of Lemma \ref{lemmaone}.

For $t\leq -1$, we have
$p_f'(t) =t^{2j-1}((2j+1) t + 2j)-f >0$,
so $p_f$ is strictly increasing on
$(-\infty,-1)$ and since $p_f(-1)= 2f <0$ this
implies that $p_f$ does not have a root in
$(-\infty,-1]$.
Since $p_f(1)= 2>0$, there exists a real
root of $p_f$ in $(-1,1)$ and all real roots of $p_f$ are contained in
this interval.


For $f=0$, $z_0=-1$ is a single root of $p_f$.
Hence there exist
numbers $\epsilon>0$, $\delta>0$ and a
neighbourhood $U(-1)=\{z\in C:\;|z+1|\leq \epsilon\}$
such that for all $f\in (-\delta,0)$
there exists exactly one root of $p_f$ in $U(-1)$.
Since the complex roots appear in conjugate pairs,
this must be a real root, hence it is
in $(-1,1)$.

In particular, its absolute value is strictly less than one.

The other $2j$ roots of $p_f$ are all outside of $U(-1)$.
If $|f|$ is sufficiently small, for all $z\not \in U(-1)$
we have
$|z-(-1)|> \epsilon > \frac{2 |f|}{ 1 - |f|}$ and
thus Lemma \ref{lemma0p} implies $|z|<1$
which finishes the proof of
Lemma \ref{lemmaone}.

The following Lemma implies a necessary condition
that must hold if all roots of $p_f$
have a  modulus that is strictly less than one.
\begin{lemma}
\label{lemma1p}
If $f<0$ and $p_f(-f)\geq 0$, there exits a real root
of $p_f$ with absolute value greater than or equal to one.
\end{lemma}
{\bf Proof.}
We have $p_f(0)= f <0$. Since $p(-f)\geq 0$,
there exists  a real root $t_\ast\in (0,-f]$
with $p_f(t_\ast=0$.
Hence $\frac{-f}{t_\ast}\geq 1$.
Let $z_1$,...,$z_{2j}$ denote the other roots.
Then we have
$|z_1 z_2 \cdots z_{2j}| = \frac{-f}{t_\ast}\geq 1$.
Thus the assertion follows.

\begin{lemma}
\label{lemma2p}
If $f<0$ and $\lambda \in ( 3 - 2\sqrt{2}, 3 + \sqrt{2})$, we have
$$p_f((\lambda |f|)^{\frac{1}{2j+1}}) >0.$$
\end{lemma}
{\bf Proof.}
We have
\begin{eqnarray*}
p_f((\lambda |f|)^{\frac{1}{2j+1}})  & = &
\lambda |f| + \lambda^\frac{2j}{2j+1} |f|^\frac{2j}{2j+1}
    - \lambda^{\frac{1}{2j+1}}  f  |f|^{\frac{1}{2j+1}} + f
\\
& = &
(\lambda - 1) |f| + \lambda^\frac{2j}{2j+1} |f|^\frac{2j}{2j+1} +
 \lambda^{\frac{1}{2j+1}} |f|^{\frac{2j + 2}{2j+1}}
\\
& = &
|f|^\frac{2j}{2j+1}
\left[
(\lambda - 1) |f|^\frac{1}{2j+1}
+\lambda^\frac{2j}{2j+1}
+
\lambda^{\frac{1}{2j+1}}
 |f|^\frac{2}{2j+1}
\right]
\\
& = &
|f|^\frac{2j}{2j+1} \;h( |f|^\frac{1}{2j+1}),
\end{eqnarray*}
where
$$h(t)= \lambda^{\frac{1}{2j+1}} t^2 + (\lambda - 1) t
+ \lambda^\frac{2j}{2j+1}.$$
Since $$\Delta =
(\lambda - 1)^2 - 4 \lambda = \lambda^2 - 6\lambda + 1
<0$$
due to our choice of $\lambda$, we have
$h(t)>0$ for all $t\in (-\infty,\infty)$ and the assertion follows.

\begin{lemma}
\label{lemma3p}
If $f<0$ and
\begin{equation}
\label{bedingung} f> - \frac{1}{2j}
\left(\frac{2j-1}{2j+1}\right)^{2j},
\end{equation}
 we have
\begin{equation}
\label{gleieins}
 p_f(- \frac{2j-1}{2j+1}) = \frac{1}{2j+1} \left[ 4j
f + 2 \left( \frac{2j-1}{2j+1} \right)^{2j}\right]
>0,
\end{equation}
\begin{equation}
\label{gleizwei}
 p_f(-(2j|f|)^\frac{1}{2j}) >0,
\end{equation}
\begin{equation}
\label{gleivier}
 p_f(-|f|^\frac{1}{2j}) \leq 0,
\end{equation}
\begin{equation}
\label{gleidrei} p_f((\frac{1}{2j}|f|)^\frac{1}{2j}) < 0
\end{equation}
and there exist three real roots of $p_f$, one of them in the
interval $(-1, - \frac{2j-1}{2j+1})$, another in the interval $( (
-(2j|f|)^\frac{1}{2j} ,\,-|f|^\frac{1}{2j} )$ and the third in $(
(\frac{1}{2j}|f|)^\frac{1}{2j} ,1)$.

For $f=- \frac{1}{2j} \left(\frac{2j-1}{2j+1}\right)^{2j}$, we have
$p_f(- \frac{2j-1}{2j+1}) =0$.
\end{lemma}
{\bf Proof.}
Let $t=- \frac{2j-1}{2j+1}$. We have
\begin{eqnarray*}
p_f(t) & = & t^{2j}(1+t) + f(1-t)
\\
& = & f( 1 +  \frac{2j-1}{2j+1}) + t^{2j} ( 1 -  \frac{2j-1}{2j+1})\\
& = &
f \frac{4j}{2j+1} + t^{2j} \frac{2}{2j+1}\\
& = & \frac{1}{2j+1} \left[ 4j  f + 2 t^{2j}\right]
\\
& > & 0
\end{eqnarray*}
which implies (\ref{gleieins}).

For $\lambda\in (1,\infty)$, we have $p_f(-(\lambda
|f|)^{\frac{1}{2j}}) = |f| [-(\lambda +1) \lambda^{\frac{1}{2j}}
|f|^{\frac{1}{2j}} + \lambda -1 ] >0$ if
$$|f| < \left(\frac{\lambda - 1}{\lambda + 1}\right)^{\frac{1}{2j}}
\frac{1}{\lambda}$$ and (\ref{gleizwei}) follows with the choice
$\lambda = 2j$.

We have $p_f(-|f|^{\frac{1}{2j}}) = |f| [-2 |f|^{\frac{1}{2j}} ] \leq 0$  and (\ref{gleivier}) follows.

For $\lambda\in (0,1)$, we have $p_f((\lambda |f|)^{\frac{1}{2j}}) =
|f| [(\lambda +1) \lambda^{\frac{1}{2j}} |f|^{\frac{1}{2j}} +
\lambda -1 ] <0$ if
$$|f| < \left(\frac{1 - \lambda}{1 + \lambda}\right)^{\frac{1}{2j}}
\frac{1}{\lambda}$$ and (\ref{gleidrei}) follows with the choice
$\lambda = \frac{1}{2j}$.

Since $p_f(-1)= 2f <0$ and $p_f(0)=f<0$ and $p_f(1)=2>0$
the assertion follows.

Lemma \ref{lemma1p} implies that we only need to consider
values of $f<0$ with $p_f(-f)<0$.

For $j=1$ this yields the sharper result given in Lemma \ref{lemmatwo}.

\begin{lemma}
\label{lemmatwo}
Let $j=1$. Then
for all $f\in (1-\sqrt{2},0)$ we have
$p_f(-f)<0$ and $p_f((-f)/5)^{1/3})>0$.
Hence there exists a root of $p_f$ in the interval
$(-f, (-f/5)^{1/3})$.
The other two roots we have a modulus that is strictly less than one.
\end{lemma}
{\bf Proof.}
We have $p(-f)=-f^3 + 2f^2 + f = -f[(f-1)^2 -2]$.
Since $f\in (1-\sqrt{2},0)$ we have $1 < (1-f)^2 < 2$
which implies $p_f(-f)<0$.
Lemma \ref{lemma2p} with $\lambda = 1/5$ implies
$p_f((-f)/5)^{1/3})>0$.
Hence there exists a root $t_\ast$ of $p_f$ in the interval
$(-f, (-f/5)^{1/3})$.

If the other roots are complex conjugate,
we call them $z$ and $\bar z$ and have
$|z|^2 = z\bar z = -f/t_\ast <1$.

Now we consider the case that the other roots are real.
Note that for $t>0$, we have $p'(t)>-f >0$
so there exists nor root that is greater than $t_\ast$.
On the other hand, for $t<-1$
we have
$p'(t)  > 1 - f >0$.
Since $p(-1)= 2 f <0$, this implies that
there is no root in $(-\infty,-1]$
hence also in this case the absolute value of all three roots is strictly less than one.
Hence the assertion follows.

For the case $j=2$ where $p_f$ is a polynomial of degree five
we only have the  result given in Lemma \ref{lemmathree}.

\begin{lemma}
\label{lemmathree}
Let $j=2$. Then
for all $f\in (-  81/2500 ,0)$ we have
$p_f(-f)<0$ and
the roots of $p_f$  have a modulus that is strictly less than one.
\end{lemma}
{\bf Proof.}
Case 1:
Suppose that $p_f$ had five real roots.
Then they would all be in the interval $(-1,1)$

Case 2:
Now we consider the case that $p_f$ has
two complex conjugate roots $z$ and $\bar z$.
Since $f \in (- 81/(4 * 625),\,0)$,
Lemma \ref{lemma3p}
implies that we have three real roots
$t_1$, $t_2$, $t_3$ such that $-1 < t_1< - \frac{3}{5}
< t_2 < - |f|^{1/4} < \frac{1}{4^{1/4}}  |f|^{1/4} < t_3$.

We have
$$|z|^2 = z\bar z = -\frac{f}{t_1\,t_2\,t_3}
<
\frac{5}{3}\frac{|f|}{\frac{1}{\sqrt{2}} |f|^{1/2}}
= \frac{5 \sqrt{2}}{3} \sqrt{|f|} <1.
$$
Hence all roots have a modulus that is strictly less than one
and the first part of
the assertion follows with Lemma \ref{lemma1p}.

The following Lemma gives a construction of values of $f$
for which a pair of complex conjugate roots of $p_f$ for $j=2$ is known.
If $|f|$ is sufficiently small and $f<0$, the remaining  three roots are all real,
so they can be easily approximated to arbitrary precision.

\begin{lemma}
\label{lemmafour}
Let $j=2$.
Let $a\geq 0$ be given.
Define the numbers
$q=\frac{a+4a^2+2a^3}{1+2a}$,
$R=q+\sqrt{q^2-4a^3}$
and
$$f=(8a^3+4a^2)R-(1+4a)R^2.$$
Let
$b=1+2a$,
$c=2ab-R$ and
$d=2 a c - Rb$.
Then we have
$$p_f(z)=(z^2 - 2a z + R)\;(z^3 + b z^2 + c z + d).$$
In particular, $p_f$ has the roots
$z_1= a + \sqrt{R-a^2}i$ and $z_2= a - \sqrt{R-a^2}i$.

\end{lemma}

{\bf Proof.}
We have
$$(z-z_1)(z-z_2)=z^2 -2a z +R.$$
Hence
\begin{eqnarray*}
(z-z_1)(z-z_2)\;(z^3 + b z^2 + c z + d) & = & (z^2 - 2a z + R)\;(z^3 + b z^2 + c z + d)
\end{eqnarray*}
$$
 =
z^5 + (b-2a)z^4 + (c-2ab + R)z^3 + (d - 2ac + Rb) z^2 + (Rc - 2ad)z + Rd
$$
\begin{eqnarray*}
& = & z^5 + z^4 +  (Rc - 2ad)z + Rd.
\end{eqnarray*}
Using the definition of $f$ we obtain the equation
$$Rd = R[4a^2(1+2a)
- R(1+4a)] = f.$$

 From the definition of $R$ we have $R^2 - 2q R+  4a^3=0$.
Hence
$$0 = (1+2a)R^2 - 2(2a^3 + 4a^2 + a)R + 4a^3(1+2a).$$
This is equivalent to the equation
$$0 = (R-2a)d + Rc.$$
Hence we have $-2ad + Rc = -Rd = -f$, thus
$$(z^2 - 2a z + R)\;(z^3 +
b z^2 + c z + d) = z^5 + z^4 -f z + f=p_f(z)$$
and the assertion follows.


\section{Exponential stability of system ${\bf S_1}$ with piecewise constant delay}
\label{piece}
We define the energy
\begin{equation}
\label{energydefinition}
E(t)=
\frac{1}{2}
\int_0^L
\left(\partial_x v(t,x)\right)^2
+
\frac{1}{c^2}
\left( \partial_t v(t,x) \right)^2 \,dx
\end{equation}
and the
energy $E_1$  by the equation
\begin{equation}
\label{energydefinition2}
E_1(t) = \sum_{j=0}^4 E\left(t + 2 j \frac{L}{c}\right).
\end{equation}
Note that $E(t)\leq E_1(t)$.

To show the exponential stability of ${\bf S_1}$,
we use the following result:

\begin{lemma}
\label{lemma2010}
Let $\lambda >0$ and the function
$E: [0,\infty)\rightarrow [0,\infty)$ be given.
Then the following two statements are equivalent:

\begin{enumerate}
\item
\label{(1)}
$E$ decays exponentially in the sense that
there exist real numbers $C_1,\,\mu \in (0,\infty)$ such that
$$ E(t) \leq C_1 \;E(0)\,\exp(-\mu t)$$
for all $t \in [0,\,\infty)$.

\item
\label{(2)}
There exist real numbers $C_2>0$ and $f \in (0,1)$ such that
the inequality
$$E(t + \,j \lambda) \leq f^j\, C_2\; E(0)$$
holds for all $t\in [0,  \lambda)$
and for all $j\in \{0,1,2,...\}$.

\end{enumerate}

\end{lemma}

{\bf Proof.}
First we show that \ref{(1)}. implies \ref{(2)}.
Assume that  \ref{(1)}. holds.
Then for all $t\in [0,\lambda)$ and all
$j\in \{0,1,2,...\}$
we have the inequality
\begin{eqnarray*}
E(t + j \lambda) & \leq & C_1 \; E(0)\; \exp  (- \mu (t + j\lambda) ) \\
& = &  C_1 \; E(0)\;\, \exp( - \mu t) \,\exp(- \mu \lambda j) \\
& \leq & C_1 \; E(0)\;\,  \exp(- \mu \lambda j) \\
& = & C_1 \; E(0)\;\, \exp ( -\lambda \mu)^j \\
& = &  \, f^j\, C_2\; E(0)\;
\end{eqnarray*}
with $C_2 = C_1$
 and $ f = \exp(-  \lambda \, \mu)$.

Now we show that  \ref{(2)}. implies \ref{(1)}.
Assume that  \ref{(2)}. holds.
For $j\in \{0,1,2,...\}$ define $t_j = j \lambda$.
For all $t\in [\lambda, \infty)$ there exists
$j\in \{0,1,2,...\}$ such that
$t\in [t_j, t_{j+1})$.
Hence we can write
$t= t_j + s $, with $s\in [0,  \lambda)$.

Define
$$\mu = - \frac{{\rm ln}(f)}{\lambda}.$$
Then ${\rm ln}(f) = -  \lambda \, \mu$.
Let $C_1 = C_2 \exp(\lambda\,\mu)$.
Then  \ref{(2)}. implies the inequality
\begin{eqnarray*}
E(t)   =  E(s + t_j) & \leq & f^j\, C_2 \; E(0)\;\\
& = & \exp(j{\rm ln}(f))\, C_2 \; E(0)\;\\
& = & \exp(- j  \lambda \mu)\, C_1 \exp(- \lambda \, \mu)\; E(0)\;\\
& = &
C_1 \exp( - \mu t_j) \, \exp (  -  \mu  \lambda) \; E(0)\;\\
& \leq & C_1 \; E(0)\;\, \exp(- \mu(t_j + s))
\\
& =  & C_1 \; E(0)\;\, \exp(-\mu \,t)
\end{eqnarray*}
and the assertion follows.

\begin{theorem}
\label{theorem0}
Let 
%
$$f_0 = \frac{-519801-761 \sqrt{467857}}{303170688}=0.00343...$$
Assume that the delay $\delta$ is piecewise constant and that for all $t\geq 0$ we have
$\delta(t)\in  \{4L/c,\; 8L/c\}$.

Then there exists a neightbourhood $U$ of $f_0$ such that for all $f\in U$
System ${\bf S_1}$ with $\iota=2$ is exponentially stable in the sense that
the energy decays exponentially.
In fact
there exists a  constant $C_0>0$
that is independend of the initial state $(y_0,y_1)$ and
a constant $L<1$
such that for all
$j\in \{0,1,2,...\}$
and for all $t\in [0, 2 L/c)$ we have the inequality
$$
E_1( t + 2j \frac{L}{c})
\leq
L^j
C_0 \;E_1(0).
$$
\end{theorem}
{\bf Proof.}
Theorem \ref{solution} states that system ${\bf S_1}$
has a solution
for which we can compute the corresponding energy
defined in (\ref{energydefinition}) as
\begin{eqnarray*}
E(t)
& = &
\int_0^L
\alpha'(x+ct )^2
 +  \alpha'( -x + ct )^2
\,dx
\\
& = &
\int_{-L}^L \alpha'(x + ct)^2\, dx.
\end{eqnarray*}
Let $h=2L$.
Let $\lambda_i$, $i\in \{1,2,3,4,5\}$ denote the eigenvalues
of the matrix $B_2$ from system (\ref{21}).
Assume that we have
$|\lambda_5| \leq |\lambda_4|\leq |\lambda_3|\leq |\lambda_2|\leq |\lambda_1|$.

Note that for
$a=1/36$,
$q=\frac{a+4a^2+2a^3}{1+2a}$ and
$R=q+\sqrt{q^2-4a^3}$ we have
$f_0=4(2a^3+a^2)R-(1+4a)R^2$. Hence for $f=f_0$,
due to Lemma \ref{lemmafour} we have the eigenvalues
$ a \pm \sqrt{R-a^2}i$.
Due to Lemma  \ref{lemma3p}
the other three eigenvalues are real and can be approximated as the
roots of the polynomial of degree that is given in Lemma \ref{lemmafour}, namely
$$z^3 + \frac{19}{18} z^2 + \frac{723-\sqrt{467857}}{24624}z -  \frac{3244 + 5 \sqrt{467857}}{110808}.$$

Define the corresponding eigenvectors
$$
s_i=
\frac{1}{\sqrt{1+\lambda_i^2+ \lambda_i^4+\lambda_i^6+ \lambda_i^8}}
\left(
\begin{array}{r}
\lambda_i^4\\
\lambda_i^3\\
\lambda_i^2\\
\lambda_i\\
1
\end{array}
\right).
$$
and the matrix
$$
V_2 =
\left(
\begin{array}{r|r|r|r|r}
s_1 & s_2 & s_3 & s_4 & s_5
\end{array}
\right)
$$
Choose the functions
$c_1(s)$, $c_2(s)$, $c_3(s)$, $c_4(s)$, $c_5(s)$ such that for
$s\in (-L,9L)$ almost everywhere we have
$$
\left(
\begin{array}{r}
\alpha'(s+4h)
\\
\alpha'(s+3h)
\\
\alpha'(s+2h)
\\
\alpha'(s+h)
\\
\alpha'(s)
\end{array}
\right)
=
V_2
\;
\left(
\begin{array}{r}
c_1(s)
\\
c_2(s)
\\
c_3(s)
\\
c_4(s)
\\
c_5(s)
\end{array}
\right).
$$
Since the matrix is invertible and
$\alpha'\in L^2_{loc}(-L,\infty)$ this implies
$c_1$, $c_2$, $c_3$, $c_4$, $c_5$ in $ L^2_{loc}(-L,9L)$.
The functions $c_i$ are the coefficients
of the representation as a linear combination of the
eigenvectors of the matrix $B_2$.
%
Then for all natural numbers $j\in \{0,1,2,...\}$
due to
(\ref{21}) and (\ref{20}) we have the representation
$$
\left(
\begin{array}{r}
\alpha'(s+4h+jh)
\\
\alpha'(s+3h+jh)
\\
\alpha'(s+2h+jh)
\\
\alpha'(s+h+jh)
\\
\alpha'(s+jh)
\end{array}
\right)
=
\sum_{i=1}^5 \gamma_{i,j}(s) c_i(s) s_i,
$$
where
$$\left(
\begin{array}{r}
\gamma_{1,j+1}(s)
\\
\gamma_{2,j+1}(s)
\\
\gamma_{3,j+1}(s)
\\
\gamma_{4,j+1}(s)
\\
\gamma_{5,j+1}(s)
\end{array}
\right)
=M(s)
\left(
\begin{array}{r}
\gamma_{1,j}(s)
\\
\gamma_{2,j}(s)
\\
\gamma_{3,j}(s)
\\
\gamma_{4,j}(s)
\\
\gamma_{5,j}(s)
\end{array}
\right)
$$
with the matrix
$$M(s)=\left\{
\begin{array}{rrr}
V_2^{-1} B_2 V_2 & {\rm if }& \delta(s/c)= 8L/c, \\
V_2^{-1} B_1 V_2 & {\rm if }& \delta(s/c)= 4L/c.
\end{array}
\right.
$$
By our construction, the matrix
$V_2^{-1} B_2 V_2=D_2 $ is a diagonal matrix
that contains the numbers $\lambda_i$ as
diagonal elements.
Due to Lemma  \ref{lemmathree} this implies that
$\|D_2\|_1<1$. (In fact, we have $\|D_2\|_1<0.994$.)

Let $H_1=V_2^{-1} B_1 V_2$.
Then numerical computations show that
$\|H_1\|_1<1$. (In fact, we have $\|H_1\|_1<0.997$.)

Define $L_0=\max\{\|D_2\|_1,\, \|H_1\|_1\}<1$.
Then we have  the inequality
\begin{equation}
\label{ungleichung0}
\left\|
\left(
\begin{array}{r}
\gamma_{1,j}(s)
\\
\gamma_{2,j}(s)
\\
\gamma_{3,j}(s)
\\
\gamma_{4,j}(s)
\\
\gamma_{5,j}(s)
\end{array}
\right)
\right\|_1
\leq
L_0^j
\left\|
\left(
\begin{array}{r}
\gamma_{1,0}(s)
\\
\gamma_{2,0}(s)
\\
\gamma_{3,0}(s)
\\
\gamma_{4,0}(s)
\\
\gamma_{5,0}(s)
\end{array}
\right)
\right\|_1
=5 L_0^j.
\end{equation}
This implies the inequality
\begin{eqnarray*}
\left\|
\left(
\begin{array}{r}
\alpha'(s+4h+jh)
\\
\alpha'(s+3h+jh)
\\
\alpha'(s+2h+jh)
\\
\alpha'(s+h+jh)
\\
\alpha'(s+jh)
\end{array}
\right)
\right\|_2
&\leq &
\|\sum_{i=1}^5 \gamma_{i,j}(s) c_i(s) s_i\|_2
=
\left\| V_2\left(
\begin{array}{r}
\gamma_{1,j}(s)\,c_1(s)
\\
\gamma_{2,j}(s)\,c_2(s)
\\
\gamma_{3,j}(s)\,c_3(s)
\\
\gamma_{4,j}(s)\,c_4(s)
\\
\gamma_{5,j}(s)\,c_5(s)
\end{array}
\right)
 \right\|_2
\\
& \leq &
\|V_2\|_2
\left\|
\left(
\begin{array}{r}
\gamma_{1,j}(s)\,c_1(s)
\\
\gamma_{2,j}(s)\,c_2(s)
\\
\gamma_{3,j}(s)\,c_3(s)
\\
\gamma_{4,j}(s)\,c_4(s)
\\
\gamma_{5,j}(s)\,c_5(s)
\end{array}
\right)
\right\|_2
\\
& \leq &
\|V_2\|_2
\left\|
\left(
\begin{array}{r}
c_1(s)
\\
c_2(s)
\\
c_3(s)
\\
c_4(s)
\\
c_5(s)
\end{array}
\right)
\right\|_2^{1/2}
\,
\left\|
\left(
\begin{array}{r}
\gamma_{1,j}(s)
\\
\gamma_{2,j}(s)
\\
\gamma_{3,j}(s)
\\
\gamma_{4,j}(s)
\\
\gamma_{5,j}(s)
\end{array}
\right)
\right\|_2^{1/2}
\\
& \leq &
\|V_2\|_2
\left\|
\left(
\begin{array}{r}
c_1(s)
\\
c_2(s)
\\
c_3(s)
\\
c_4(s)
\\
c_5(s)
\end{array}
\right)
\right\|_2^{1/2}
\left\|
\left(
\begin{array}{r}
\gamma_{1,j}(s)
\\
\gamma_{2,j}(s)
\\
\gamma_{3,j}(s)
\\
\gamma_{4,j}(s)
\\
\gamma_{5,j}(s)
\end{array}
\right)
\right\|_1^{1/2}
\\
& \leq &
\sqrt{5}
\|V_2\|_2
\,L_0^{j/2}
\left\|
\left(
\begin{array}{r}
c_1(s)
\\
c_2(s)
\\
c_3(s)
\\
c_4(s)
\\
c_5(s)
\end{array}
\right)
\right\|_2^{1/2}.
\end{eqnarray*}

Let $t\in [0, 10L/c)$.
For the energy $E_1$ we have the equation
\begin{eqnarray*}
E_1(t + 2j \frac{L}{c})
& = &
\int_{-L}^{9L} \alpha'(x +  c t + j h)^2\, dx
\\
& \leq &
\int_{-L}^{L}
\left\|
\left(
\begin{array}{r}
\alpha'(s +ct +4h+jh)
\\
\alpha'(s + ct +3h+jh)
\\
\alpha'(s + ct +2h+jh)
\\
\alpha'(s + ct +h+jh)
\\
\alpha'(s + ct +jh)
\end{array}
\right)
\right\|_2^2
\,ds
\\
& \leq &
{5}
\|V_2\|_2^2
\,L_0^{j}
\int_{-L}^{L}
\left\|
\left(
\begin{array}{r}
c_1(s)
\\
c_2(s)
\\
c_3(s)
\\
c_4(s)
\\
c_5(s)
\end{array}
\right)
\right\|_2
\,ds
\\
& \leq &
L^j_0 C_0 \;E_1(0)
\end{eqnarray*}
which implies the exponential decay for $f=f_0$ due to Lemma \ref{lemma2010}.
Due to continuity, we find a neighbourhood $U$ of $f_0$ such that
for all $f\in U$ we
have $\|D_2(f)\|_1<1$ and $\|H_1(f)\|_1<1$ and this yields the assertion.

\section{Exponential stability of system ${\bf S_1}$ with constant delay}

\begin{theorem}
\label{theorem1}
For all $\iota\in \{0,1,2,...\}$ there exists a number $\delta_\iota>0$ such that for all
$f\in (-\delta_\iota,0)$
System ${\bf S_1}$ 
with the constant delay $\delta(t)= 4 \iota L/c$ 
is exponentially stable in the sense that
the energy decays exponentially 
In fact
there exists a constant $C_0>0$
that only depends on the initial state $(y_0,y_1)$ and $f$
such that for all
$j\in \{0,1,2,...\}$
and for all $t\in [0, 2 L/c)$ we have the inequality
$$
E( t + 2j \frac{L}{c})
\leq
f^j
C_0 \;E(0).
$$
\end{theorem}
\begin{remark}
Note that for the
corresponding feedback law without delay
\begin{equation}
\label{without}
c v_x(t,L)  =  f \;  v_t\left(t,L\right),\; t>0
\end{equation}
with $f=-1$, the energy is controlled to
zero in finite time.
\end{remark}

\subsection{Proof of Theorem \ref{theorem1}} 
Let $\iota\in \{0,1,2,...\}$ be given.
Define  the
characteristic
polynomial
$p_f(t)$ as in Section \ref{charac} with $j=\iota$.
Lemma \ref{lemmaone} states that
there exists a number $\delta_\iota>0$, such that for all
$f\in (-\delta_\iota,0)$, all roots of $p_f$ have
a modulus that is strictly less than one.
The proof uses the fact that from (\ref{2.8}) we get an explicit
representation of $\alpha'$.
Let $z_1$,....,$z_{2\iota+1}$ denote the roots of $p_f$.

Theorem \ref{solution} states that system ${\bf S_1}$
has a solution
for which we can compute the corresponding energy
defined in (\ref{energydefinition}) as
\begin{eqnarray*}
E(t)=\int_{-L}^L \alpha'(x + ct)^2\, dx.
\end{eqnarray*}
Let $h=2L$. For $x \geq (4\iota +1)L$
equation (\ref{2.8})
yields the equation
\begin{equation}
\alpha'(s) +  \alpha'(s - h) - f \alpha'(s - 2 \iota h)
+ f  \alpha'(s - (2 \iota +1) h)
\label{ld1}
=0
.
\end{equation}

Using the usual method for linear difference equations
we obtain an explicit representation of the solution
 $\alpha'\in L^2_{loc}(-L,\infty)$
of (\ref{ld1}).
Choose the functions
$c_1(s)$, $c_2(s)$,...,$c_{2\iota+1}(s)$ such that for  $s\in (-L,(2\iota+1)L)$ almost everywhere we have
$$
\left(
\begin{array}{r}
\alpha'(s)
\\
\alpha'(s+h)
\\
\hdots
\\
\alpha'(s+2\iota h)
\end{array}
\right)
=
\left(
\begin{array}{rrrr}
1 & 1 &  \hdots & 1 
\\
z_1 & z_2 & \hdots & z_{2\iota + 1}\\
\\
\vdots & & & \vdots \\
z_1^{2\iota}& z_2^{2\iota} & \hdots & z_{2\iota + 1}^{2\iota}
\end{array}
\right)
\;
\left(
\begin{array}{r}
c_1(s)
\\
c_2(s)
\\
\vdots
\\
c_{2\iota + 1}(s)
\end{array}
\right).
$$
Since the matrix is invertible and
$\alpha'\in L^2_{loc}(-L,\infty)$ this implies
$c_k$ in $ L^2_{loc}(-L,L)$.
%
Then for all natural numbers $j\in 1,2,...
$
 we have the representation
$$\alpha'(s + jh) = c_1(s) z_1^j + c_2(s) z_2^j +
... + c_{2\iota+1}(s) \, z_{2\iota+1}^j$$

Let $M_0 = \max\{|z_1|,|z_2|,...,|z_{2\iota+1}|$.
Then we have  the inequality
\begin{equation}
\label{ungleichung}
|\alpha'(s + jh)|\leq M_0^j (|c_1(s)| + |c_2(s)|
+ ... +  |c_{2\iota+1}(s)|).
\end{equation}


Let $t\in [0, 2L/c)$.
For the energy  we obtain the inequality
\begin{eqnarray*}
E(t + 2j \frac{L}{c})
& = &
\int_{-L}^L \alpha'(x + jh + c t)^2\, dx
\\
& = & \int_{-L}^{L-ct}
\alpha'(x + ct +  jh )^2\, dx
\\
& + &
 \int_{L-ct}^L
\alpha'(x + ct - 2L +  (j+1)h )^2\, dx
\\
& = & \int_{-L}^{L-ct}
\left|c_1(x + ct)z_1^j +
...+c_{2\iota+1}(x + ct)z_{2\iota+1}^j
\right|^2\, dx
\\
& + &
 \int_{L-ct}^L
\left|c_1(x + ct- 2L)z_1^{j+1} + 
...+c_{2\iota+1}(x + ct-2L)z_{2\iota+1}^{j+1}
\right|^2\, dx
\\
& \leq &
2 \,M_0^j
\,\int_{-L}^L\left(|c_1(s)| + |c_2(s)|+
...
+|c_{2\iota+1}(s)|
\right)^2\, ds
\end{eqnarray*}
which implies the exponential decay due to Lemma \ref{lemma2010}.







\section{Conclusion}
In this paper we have  considered feedback laws
that use observations from the
past for the boundary control  of the
considered systems
that are governed by the wave equation.
In this way, there is enough time for the
processing of the feedback law in practice.

We have shown that if the feedback parameters
are chosen appropriately the
feedback laws with constant delay
lead to exponential decay of the energy of the
vibrating systems if the delay is an integer
multiple of $4L/c$.

Moreover, we have shown that if the delay is piecewise
constant with values in $4L/c$, $8L/c$, the system
also decays exponentially if the feedback parameter is
chosen appropriately.



%
%
%

\noindent
 {\bf Acknowledgement}
This paper was supported by the PROCOPE program of DAAD, D/0811409.
This paper took benefit from discussions during the
meeting
{\em Partial differential equations, optimal design and numerics 2009}
at the BENASQUE Center for Science Pedro Pascual.
%
%
%

%
%


\end{document}